\documentclass[10pt]{amsart}
\usepackage{amssymb,amsmath,amsfonts,latexsym,amsthm,geometry,graphicx}
\usepackage{amsmath}
\usepackage{amsfonts}
\usepackage{amssymb}
\usepackage{graphicx}
\providecommand{\U}[1]{\protect\rule{.1in}{.1in}}
\providecommand{\U}[1]{\protect\rule{.1in}{.1in}}
\providecommand{\U}[1]{\protect\rule{.1in}{.1in}}
\newtheorem{theorem}{Theorem}[section]

\theoremstyle{definition}

\begin{document}
\title[Optimal constants for a mixed Littlewood type inequality]{Optimal constants for a mixed Littlewood type inequality}
\author[T. Nogueira]{Tony Nogueira}
\address[T. Nogueira]{Departamento de Matem\'{a}tica\\
\indent Universidade Federal da Para\'{\i}ba\\
\indent 58.051-900 - Jo\~{a}o Pessoa, Brazil.}
\email{tonykleverson@gmail.com}
\author[D. N\'{u}\~{n}ez]{Daniel N\'{u}\~{n}ez-Alarc\'{o}n}
\address[D. N\'{u}\~{n}ez-Alarc\'{o}n]{Departamento de Matem\'{a}tica\\
\indent Universidade Federal de Pernambuco\\
\indent50.740-560 - Recife, Brazil\\
\indent\& Department of Mathematical Sciences\\
\indent Kent State University\\
\indent Kent, Ohio 44242, USA}
\email{danielnunezal@gmail.com}
\author[D. Pellegrino]{Daniel Pellegrino}
\address[D. Pellegrino]{Departamento de Matem\'{a}tica \\
Universidade Federal da Para\'{\i}ba \\
58.051-900 - Jo\~{a}o Pessoa, Brazil.}
\email{pellegrino@pq.cnpq.br and dmpellegrino@gmail.com}
\thanks{T. Nogueira was supported by Capes, D. N\'{u}\~{n}ez-Alarc\'{o}n was supported
by Capes, Grant 000785/2015-06, and D. Pellegrino was supported by CNPq.}
\keywords{Absolutely summing operators; Hardy--Littlewood inequality; Bohnenblust--Hille
inequality; Multiple summing operators}
\subjclass[2010]{11Y60, 46G25.}

\begin{abstract}
For $p\in\lbrack2,\infty]$ a mixed Littlewood-type inequality asserts that
there is a constant $C_{(m),p}\geq1$ such that
\[
\left(  \sum_{i_{1}=1}^{\infty}\left(  \sum_{i_{2},...,i_{m}=1}^{\infty
}|T(e_{i_{1}},...,e_{i_{m}})|^{2}\right)  ^{\frac{1}{2}\frac{p}{p-1}}\right)
^{\frac{p-1}{p}}\leq C_{(m),p}\Vert T\Vert
\]
for all continuous real-valued $m$-linear forms on $\ell_{p}\times c_{0}
\times\dots\times c_{0}$ (when $p=\infty$, $\ell_{p}$ is replaced by $c_{0})$.
We prove that for $p>2.18006$ the optimal constants $C_{(m),p}$ are $\left(
2^{\frac{1}{2}-\frac{1}{p}}\right)  ^{m-1}.$ When $p=\infty,$ we recover the
best constants of the mixed $\left(  \ell_{1},\ell_{2}\right)  $-Littlewood inequality.

\end{abstract}
\maketitle

\section{Introduction}

The Hardy--Littlewood inequality (\cite{hardy}, 1934) is a continuation of
famous works of Littlewood (\cite{LLL}, 1930) and Bohnenblust and Hille
(\cite{bh}, 1931) and can be stated as follows:

\begin{itemize}
\item \cite[Theorems 2 and 4]{hardy} If $p,q\geq2$ are such that
\[
\frac{1}{2}<\frac{1}{p}+\frac{1}{q}<1
\]
then there is a constant $C_{p,q}\geq1$ such that
\begin{equation}
\left(  \sum\limits_{j,k=1}^{\infty}\left\vert A(e_{j},e_{k})\right\vert
^{\frac{pq}{pq-p-q}}\right)  ^{\frac{pq-q-p}{pq}}\leq C_{p,q}\left\Vert
A\right\Vert \label{yh}
\end{equation}
for all continuous bilinear forms $A:\ell_{p}\times\ell_{q}\rightarrow
\mathbb{R}$ (or $\mathbb{C}$). Moreover the exponent $\frac{pq}{pq-p-q}$ is optimal.

\item \cite[Theorems 1 and 4]{hardy} If $p,q\geq2$ are such that
\[
\frac{1}{p}+\frac{1}{q}\leq\frac{1}{2}
\]
then there is a constant $C_{p,q}\geq1$ such that
\begin{equation}
\left(  \sum\limits_{j,k=1}^{\infty}\left\vert A(e_{j},e_{k})\right\vert
^{\frac{4pq}{3pq-2p-2q}}\right)  ^{\frac{3pq-2p-2q}{4pq}}\leq C_{p,q}
\left\Vert A\right\Vert \label{yh6}
\end{equation}
for all continuous bilinear forms $A:\ell_{p}\times\ell_{q}\rightarrow
\mathbb{R}$ (or $\mathbb{C}$). Moreover the exponent $\frac{4pq}{3pq-2p-2q}$
is optimal.
\end{itemize}

Above and henceforth, as usual in this field, when $p$ and/or $q$ is infinity,
we consider $c_{0}$ instead of $\ell_{p}$ and/or $\ell_{q}.$

As mentioned in \cite[Theorem 1]{tonge} an unified version of the above two
results of Hardy and Littlewood asserts that there is a constant $C_{p,q}
\geq1$ such that
\begin{equation}
\left(  \sum\limits_{j=1}^{\infty}\left(  \sum\limits_{k=1}^{\infty}\left\vert
A(e_{j},e_{k})\right\vert ^{2}\right)  ^{\frac{\lambda}{2}}\right)  ^{\frac
{1}{\lambda}}\leq C_{p,q}\left\Vert A\right\Vert \label{dan}
\end{equation}
with $\lambda=\frac{pq}{pq-p-q},$ for all continuous bilinear forms
$A:\ell_{p}\times\ell_{q}\rightarrow$ $\mathbb{R}$ \ (in fact, in
\cite[Theorem 1]{tonge} just the complex case is considered, but for a general
approach including the real case we refer to \cite{was}; moreover the
exponents are optimal). The recent years witnessed an increasing interest in
the study of summability of multilinear operators (see, for instance,
\cite{botelho,popa,rueda}) and in estimating constants of the multilinear and
polynomial Hardy--Littlewood and related inequalities (see \cite{bayart, a22,
a1, a2, dimant, diniz, diana}). Perhaps the main motivations are potential
applications (see, for instance, \cite{monta} for applications of the
real-valued case of the estimates of the Bohnenblust--Hille inequality and
\cite{bohr, ann} for applications of the complex-valued case).

One of the most for reaching generalizations of the Hardy--Littlewood
inequality is the following theorem (see also \cite{velanga}):

\begin{theorem}
\label{aaa}(See Albuquerque, Araujo, N\'{u}\~{n}ez, Pellegrino and Rueda \cite{www}) Let $m\geq2$ be a positive integer, $1\leq k\leq m$ and
$n_{1},\dots,n_{k}\geq1$ be positive integers such that $n_{1}+\cdots
+n_{k}=m.$ If $q_{1},...,q_{k}\in\left[  \frac{1}{1-\left(  \frac{1}{p_{1}
}+\dots+\frac{1}{p_{m}}\right)  },2\right] $ and $0\leq\frac{1}{p_{1}}
+\dots+\frac{1}{p_{m}}\leq\frac{1}{2},$ then the following assertions are equivalent: \\

$(a)$ There is a constant $C_{k}=C(k,p_{1},...,p_{m},q_{1},...,q_{k})$ such
that

$$
\left( 
\sum_{i_{1}=1}^{\infty}
\left( ...
\left(  \sum_{i_{k}=1}^{\infty}|T(e_{i_{1}}^{n_1},...,e_{i_{k}}^{n_k})|^{q_k}\right)  ^{\frac{q_{k-1}}{q_k}}
... \right)^{\frac{q_1}{q_2}}
\right)^{\frac{1}{q_1}}
\leq C_{k}\Vert T\Vert
$$

for all continuous $m$-linear forms $T:\ell _{p_{1}}\times \dots \times \ell
_{p_{m}}\rightarrow \mathbb{R}$. \newline

$(b)$ The numbers $q_{1},...,q_{k}$ satisfy
$$
\frac{1}{q_{1}}+\dots+\frac{1}{q_{k}}\leq\frac{k+1}{2}-\left(  \frac{1}{p_{1}
}+\dots+\frac{1}{p_{m}}\right).
$$
\end{theorem}

Above, the notation $e_{j}^{n_{j}}$ represents the $n_{j}$-tuple
$(e_{j},...,e_{j})$. The optimal constants of the previous inequalities are
essentially unknown. Recent works have shown that in general these constants
have a sublinear growth (see \cite{bra, a2, bohr}, and references therein).
One of the few cases in which the optimal constants are known for all $m$ is
the case of mixed $\left(  \ell_{1},\ell_{2}\right)  $-Littlewood inequality
(see \cite{daniel}):

\begin{itemize}
\item The optimal constants $C_{(m),\infty}$ satisfying
\begin{equation}
\sum_{i_{1}=1}^{\infty}\left(  \sum_{i_{2},...,i_{m}=1}^{\infty}|T(e_{i_{1}
},...,e_{i_{m}})|^{2}\right)  ^{\frac{1}{2}}\leq C_{(m),\infty}\Vert
T\Vert\label{ggtt}
\end{equation}
for all continuous real $m$-linear forms $T:c_{0}\times\dots\times
c_{0}\rightarrow\mathbb{R}$ are $2^{\frac{m-1}{2}}$.
\end{itemize}

From now on $p_{0}\approx1.84742$ is the unique real number satisfying
\begin{equation}
\Gamma\left(  \frac{p_{0}+1}{2}\right)  =\frac{\sqrt{\pi}}{2}. \label{haag}
\end{equation}

Our main result provides the optimal constants of a Hardy--Littlewood-type
inequality that encompasses (\ref{ggtt}); as far as we know this is the first
time in which a Hardy--Littlewood type inequality (except for the case of
mixed $\left(  \ell_{1},\ell_{2}\right)  $-Littlewood inequality) is proved to
have optimal constants with exponential growth:

\begin{theorem}
\label{ppp} Let $m\geq2$ be a positive integer and $p\geq\frac{p_{0}}{p_{0}
-1}\approx2.18006$. The optimal constant $C_{(m),p}$ such that
\begin{equation}
\left(  \sum_{i_{1}=1}^{\infty}\left(  \sum_{i_{2},...,i_{m}=1}^{\infty
}|T(e_{i_{1}},...,e_{i_{m}})|^{2}\right)  ^{\frac{1}{2}\frac{p}{p-1}}\right)
^{\frac{p-1}{p}}\leq C_{(m),p}\Vert T\Vert, \label{777}
\end{equation}
for all continuous $m$-linear forms $T:\ell_{p}\times c_{0}\times\dots\times
c_{0}\rightarrow\mathbb{R}$ is $\left(  2^{\frac{1}{2}-\frac{1}{p}}\right)
^{m-1}$.
\end{theorem}

\begin{figure}[tbh]
\centering
\includegraphics[scale = .5]{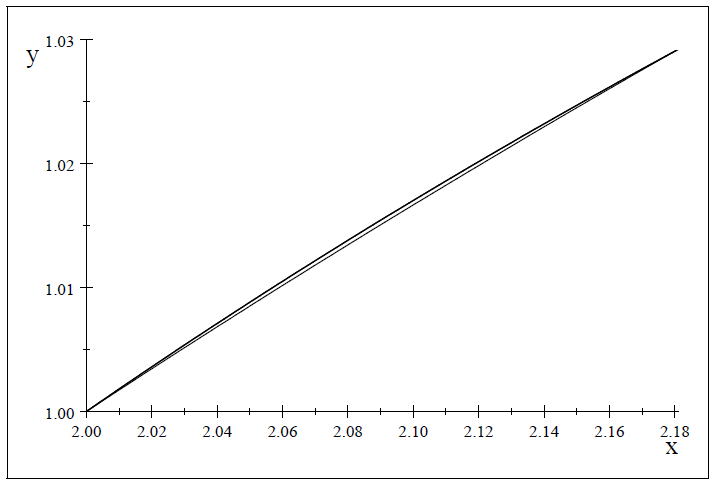}\caption{Plots of the functions
$A_{\frac{x}{x-1}}^{-1}$ and $2^{\frac{1}{2}-\frac{1}{x}}$, for $x\in
\lbrack2,\frac{p_{0}}{p_{0}-1}\rbrack$}
\label{one}
\end{figure}

\begin{figure}[tbh]
\centering
\includegraphics[scale = .5]{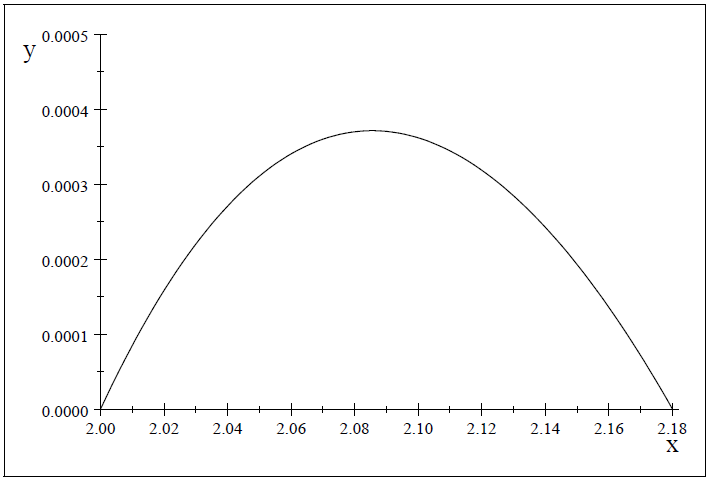}\caption{Plot of the
function $\left(  A_{\frac{x}{x-1}}^{-1}-2^{\frac{1}{2}-\frac{1}{x}}\right)
$, for $x\in\lbrack2,\frac{p_{0}}{p_{0}-1}]$}
\label{two}
\end{figure}

Note that the above Hardy--Littlewood type inequality holds for $p\geq2$ (see
Theorem \ref{aaa}). When $p=2$ it is simple to prove that the optimal
constants are $C_{(m),p}=1$. As a consequence of the arguments of our proof of
Theorem \ref{ppp} we remark that for $2<p<\frac{p_{0}}{p_{0}-1}$ the optimal
constants still have exponential growth; so an eventual decrease on the order
of the growth when $p\rightarrow2$ does not happen. Moreover, for
$2<p<\frac{p_{0}}{p_{0}-1}\approx2.18006$, the difference between the bases in
the exponential upper and lower estimates of $C_{(m),p}$ is not bigger than
$4\cdot10^{-4}$ (see the figures \ref{one} and \ref{two}).

In the final section we also provide upper and lower estimates for the sharp
constants $C_{p,\infty}$ of the real case of (\ref{yh6}), showing that
\[
2^{\frac{1}{2}-\frac{1}{p}}\leq C_{p,\infty}\leq2^{\frac{1}{2}-\frac{1}{2p}}
\]
for all $p\geq\frac{p_{0}}{p_{0}-1}\approx2.18006.$ This result recovers, in
particular, the optimality of the constant $\sqrt{2}$ of the real case of the
Littlewood's $4/3$ inequality obtained in \cite{diniz}.

\section{The proof of Theorem \ref{ppp}}

The Khinchine inequality (see \cite{Di}) asserts that, for any $0<q<\infty$,
there are positive constants $A_{q}$, $B_{q}$ such that regardless of the
scalar sequence $(a_{j})_{j=1}^{n}$ we have
\[
A_{q}\left(  \sum_{j=1}^{n}|a_{j}|^{2}\right)  ^{\frac{1}{2}}\leq\left(
\int_{0}^{1}\left\vert \sum_{j=1}^{n}a_{j}r_{j}(t)\right\vert ^{q}dt\right)
^{\frac{1}{q}}\leq B_{q}\left(  \sum_{j=1}^{n}|a_{j}|^{2}\right)  ^{\frac
{1}{2}},
\]
where $r_{j}$ are the Rademacher functions. For real scalars, U. Haagerup
\cite{Ha} proved that if $p_{0}$ is the number defined in (\ref{haag}) then
\[
A_{q}=\sqrt{2}\left(  \frac{\Gamma\left(  \frac{q+1}{2}\right)  }{\sqrt{\pi}
}\right)  ^{\frac{1}{q}},\ \ \text{ for }1.84742\approx p_{0}<q<2
\]
and
\[
A_{q}=2^{\frac{1}{2}-\frac{1}{q}},\ \ \text{ for }1\leq q\leq p_{0}
\approx1.84742.
\]
Let $T:\ell_{p}\times c_{0}\times\dots\times c_{0}\rightarrow\mathbb{R}$ be a
continuous $m$-linear form. By the Khinchine inequality for multiple sums (see
\cite{popa3}) we know that
\begin{align*}
&  \left(  \sum_{i_{1}=1}^{\infty}\left(  \sum_{i_{2},...,i_{m}=1}^{\infty
}|T(e_{i_{1}},...,e_{i_{m}})|^{2}\right)  ^{\frac{1}{2}\frac{p}{p-1}}\right)
^{\frac{p-1}{p}}\\
&  \quad\leq(A_{\frac{p}{p-1}}^{-1})^{m-1}\left(  \sum_{i_{1}=1}^{\infty}
\int_{[0,1]^{m-1}}\left\vert \sum_{i_{2},...,i_{m}}^{\infty}r_{i_{2}}
(t_{2})\cdots r_{i_{m}}(t_{m})T(e_{i_{1}},...,e_{i_{m}})\right\vert ^{\frac
{p}{p-1}}dt_{2}\cdots dt_{m}\right)  ^{\frac{p-1}{p}}\\
&  \quad=(A_{\frac{p}{p-1}}^{-1})^{m-1}\left(  \int_{[0,1]^{m-1}}\sum
_{i_{1}=1}^{\infty}\left\vert T\left(  e_{i_{1}},\sum_{i_{2}=1}^{\infty
}r_{i_{2}}(t_{2})e_{i_{2}},...,\sum_{i_{m}=1}^{\infty}r_{i_{m}}(t_{m}
)e_{i_{m}}\right)  \right\vert ^{\frac{p}{p-1}}dt_{2}\cdots dt_{m}\right)
^{\frac{p-1}{p}}\\
&  \quad\leq(A_{\frac{p}{p-1}}^{-1})^{m-1}\left(  \int_{[0,1]^{m-1}}\left\Vert
T\left(  \cdot,\sum_{i_{2}=1}^{\infty}r_{i_{2}}(t_{2})e_{i_{2}},...,\sum
_{i_{m}=1}^{\infty}r_{i_{m}}(t_{m})e_{i_{m}}\right)  \right\Vert ^{\frac
{p}{p-1}}dt_{2}\cdots dt_{m}\right)  ^{\frac{p-1}{p}}\\
&  \quad\leq(A_{\frac{p}{p-1}}^{-1})^{m-1}\sup_{t_{2},...,t_{m}\in\lbrack
0,1]}\left\Vert T\left(  \cdot,\sum_{i_{2}=1}^{\infty}r_{i_{2}}(t_{2}
)e_{i_{2}},...,\sum_{i_{m}=1}^{\infty}r_{i_{m}}(t_{m})e_{i_{m}}\right)
\right\Vert \\
&  \quad\leq(A_{\frac{p}{p-1}}^{-1})^{m-1}\Vert T\Vert=(2^{\frac{1}{2}
-\frac{1}{p}})^{m-1}\Vert T\Vert
\end{align*}
whenever $p\geq\frac{p_{0}}{p_{0}-1}\approx2.18006.$ Now let us show that
$(2^{\frac{1}{2}-\frac{1}{p}})^{m-1}$ is the best possible constant. Let
$T_{2}:\ell_{p}^{2}\times\ell_{\infty}^{2}\rightarrow\mathbb{R}$ and
$T_{2}^{x_{2}}:\ell_{p}^{2}\rightarrow\mathbb{R}$ be given by
\begin{equation}
T_{2}\left(  x_{1},x_{2}\right)  =\left(  x_{2}^{1}+x_{2}^{2}\right)
x_{1}^{1}+\left(  x_{2}^{1}-x_{2}^{2}\right)  x_{1}^{2},\label{tdos}
\end{equation}
and
\[
T_{2}^{x_{2}}\left(  x_{1}\right)  =T_{2}\left(  x_{1},x_{2}\right)  ,
\]
for each $x_{2}\in\ell_{\infty}^{2}$. Observe that
\begin{equation}
\left\Vert T_{2}\right\Vert =\sup\left\{  \left\Vert T_{2}^{x_{2}}\right\Vert
:\left\Vert x_{2}\right\Vert _{\ell_{\infty}^{2}}=1\right\}  .\label{one222}
\end{equation}
Let us estimate (\ref{one222}). Since $\left(  \ell_{p}\right)  ^{\ast}
=\ell_{\frac{p}{p-1}}$, we have
\begin{align}
\left\Vert T_{2}\right\Vert  &  =\sup\left\{  \left\Vert T_{2}^{x_{2}
}\right\Vert :\left\Vert x_{2}\right\Vert _{\ell_{\infty}^{2}}=1\right\}
\label{2211}\\
&  =\sup\left\{  \sup_{x_{1}\in B_{\ell_{p}^{2}}}\left\vert T_{2}^{x_{2}
}\left(  x_{1}\right)  \right\vert :\left\Vert x_{2}\right\Vert _{\ell
_{\infty}^{2}}=1\right\}  \nonumber\\
&  =\sup\left\{  \sup_{x_{1}\in B_{\ell_{p}^{2}}}\left\vert \left(  x_{2}
^{1}+x_{2}^{2}\right)  x_{1}^{1}+\left(  x_{2}^{1}-x_{2}^{2}\right)  x_{1}
^{2}\right\vert :\left\Vert x_{2}\right\Vert _{\ell_{\infty}^{2}}=1\right\}
\nonumber\\
&  =\sup\left\{  \left\Vert \left(  x_{2}^{1}+x_{2}^{2},x_{2}^{1}-x_{2}
^{2},0,0,...\right)  \right\Vert _{\frac{p}{p-1}}:\left\Vert x_{2}\right\Vert
_{\ell_{\infty}^{2}}=1\right\}  \nonumber\\
&  =\sup\left\{  \left(  \left\vert 1+x\right\vert ^{\frac{p}{p-1}}+\left\vert
1-x\right\vert ^{\frac{p}{p-1}}\right)  ^{\frac{1}{\frac{p}{p-1}}}:x\in
\lbrack-1,1]\right\}  =2.\nonumber
\end{align}
In order to verify the last equality, note that since
\[
\sup\left\{  \left(  \left\vert 1+x\right\vert ^{1}+\left\vert 1-x\right\vert
^{1}\right)  ^{1};x\in\lbrack-1,1]\right\}  =2,
\]
by the norm inclusion $\ell_{1}\subset\ell_{\frac{p}{p-1}}$ for $p\in
\lbrack2,\infty)\,,$ we have $\left\Vert \cdot\right\Vert _{\ell_{\frac
{p}{p-1}}}\leq\left\Vert \cdot\right\Vert _{\ell_{1}}$. Therefore, for
$p\in\lbrack2,\infty)$ we have
\[
\sup\left\{  \left(  \left\vert 1+x\right\vert ^{\frac{p}{p-1}}+\left\vert
1-x\right\vert ^{\frac{p}{p-1}}\right)  ^{\frac{p-1}{p}};x\in\lbrack
-1,1]\right\}  \leq\sup\left\{  \left(  \left\vert 1+x\right\vert
^{1}+\left\vert 1-x\right\vert ^{1}\right)  ^{1};x\in\lbrack-1,1]\right\}  =2.
\]
On the other hand, it is obvious that
\[
\sup\left\{  \left(  \left\vert 1+x\right\vert ^{\frac{p}{p-1}}+\left\vert
1-x\right\vert ^{\frac{p}{p-1}}\right)  ^{\frac{p-1}{p}};x\in\lbrack
-1,1]\right\}  \geq\left(  \left\vert 1+1\right\vert ^{\frac{p}{p-1}
}+\left\vert 1-1\right\vert ^{\frac{p}{p-1}}\right)  ^{\frac{p-1}{p}}=2.
\]

In order to show that $\left(  2^{\frac{1}{2}-\frac{1}{p}}\right)  ^{m-1}$ is
the best possible constant satisfying (\ref{777}), let $T_{2}$ be as in
(\ref{tdos}) and define for all $m\geq3$ the $m$-linear operator $T_{m}
:\ell_{p}^{2^{m-1}}\times\ell_{\infty}^{2^{m-1}}\times\dots\times\ell_{\infty
}^{2^{m-1}}\rightarrow\mathbb{R}$ by
\begin{align*}
T_{m}(x_{1},....,x_{m})=  &  (x_{m}^{1}+x_{m}^{2})T_{m-1}(x_{1},...,x_{m-1})\\
&  +(x_{m}^{1}-x_{m}^{2})T_{m-1}(S_{p}^{2^{m-2}}(x_{1}),S_{0}^{2^{m-2}}
(x_{2}),S_{0}^{2^{m-3}}(x_{3})...,S_{0}^{2}(x_{m-1})),
\end{align*}
where $x_{1}\in\ell_{p}^{2^{m-1}}$, $x_{k}\in\ell_{\infty}^{2^{m-1}}$ for all
$k=2,...,m$, and \ $S_{p}:\ell_{p}\rightarrow\ell_{p}$ and $S_{0}
:c_{0}\rightarrow c_{0}$ are the backward shifts. By induction on $m\geq2$ we
shall show that
\[
\Vert T_{m}\Vert=2^{m-1}.
\]
The case $m=2$ is already done in (\ref{2211}). Let us suppose that
$\left\Vert T_{m-1}\right\Vert =2^{(m-1)-1}$. Therefore,
\begin{align*}
|T_{m}(x_{1},\ldots,x_{m})|\leq &  |x_{m}^{1}+x_{m}^{2}||T_{m-1}(x_{1}
,\ldots,x_{m-1})|\\
&  +|x_{m}^{1}-x_{m}^{2}||T_{m-1}(S_{p}^{2^{m-2}}(x_{1}),S_{0}^{2^{m-2}}
(x_{2}),S_{0}^{2^{m-3}}(x_{3})...,S_{0}^{2}(x_{m-1}))|\\
\leq &  2^{m-2}[|x_{m}^{1}+x_{m}^{2}|\Vert x_{1}\Vert_{\ell_{p}^{2^{m-1}}
}\cdots\Vert x_{m-1}\Vert_{\ell_{\infty}^{2^{m-1}}}\\
&  +|x_{m}^{1}-x_{m}^{2}|\Vert S_{p}^{2^{m-2}}(x_{1})\Vert_{\ell_{p}^{2^{m-1}
}}\Vert S_{0}^{2^{m-2}}(x_{2})\Vert_{\ell_{\infty}^{2^{m-1}}}\Vert
S_{0}^{2^{m-3}}(x_{3})\Vert_{\ell_{\infty}^{2^{m-1}}}\cdots\Vert S_{0}
^{2}(x_{m-1})\Vert_{\ell_{\infty}^{2^{m-1}}}]\\
\leq &  2^{m-2}[|x_{m}^{1}+x_{m}^{2}|+|x_{m}^{1}-x_{m}^{2}|]\Vert x_{1}
\Vert_{\ell_{p}^{2^{m-1}}}\cdots\Vert x_{m-1}\Vert_{\ell_{\infty}^{2^{m-1}}}\\
=  &  2^{m-1}\Vert x_{1}\Vert_{\ell_{p}^{2^{m-1}}}\cdots\Vert x_{m-1}
\Vert_{\ell_{\infty}^{2^{m-1}}}\max\{|x_{m}^{1}|,|x_{m}^{2}|\}\\
\leq &  2^{m-1}\Vert x_{1}\Vert_{\ell_{p}^{2^{m-1}}}\cdots\Vert x_{m}
\Vert_{\ell_{\infty}^{2^{m-1}}}.
\end{align*}
We thus have $\Vert T_{m}\Vert\leq2^{m-1}$. Now consider $a_{m}=e_{1}+e_{2}$
and note that
\begin{align*}
\Vert T_{m}\Vert &  \geq\sup\left\{  \left\vert T_{m}\left(  x_{1}
,\dots,x_{m-1},a_{m}\right)  \right\vert :x_{1}\in B_{\ell_{p}^{2^{m-1}}
},x_{2}\in B_{\ell_{\infty}^{2^{m-1}}},...,x_{m-1}\in B_{\ell_{\infty
}^{2^{m-1}}}\right\} \\
&  =2\Vert T_{m-1}\Vert=2^{m-1}
\end{align*}
and hence $\Vert T_{m}\Vert=2^{m-1}.$

Since
\[
\frac{\left(  \sum_{i_{1}}\left(  \sum_{i_{2},...,i_{m}}|T_{m}(e_{i_{1}
,...,}e_{i_{m}})|^{2}\right)  ^{\frac{1}{2}\frac{p}{p-1}}\right)  ^{\frac
{p-1}{p}}}{\Vert T_{m}\Vert}=\left(  2^{^{\frac{1}{2}-\frac{1}{p}}}\right)
^{m-1},
\]
the proof is done.

\section{Final remarks}

The same argument used in the proof of Theorem \ref{ppp} shows that for
$2<p<\frac{p_{0}}{p_{0}-1}\approx2.18006$ the optimal constants also have
exponential growth; curiously, for $p=2$ the situation is quite different and
the optimal constants are $1.$ In fact, note that the second part of the proof
(the optimality proof) holds for all $p\geq2$. Moreover, the first part of the
proof gives us the estimate $C_{(m),p}\leq\left(  A_{\frac{p}{p-1}}
^{-1}\right)  ^{m-1}.$ We thus have, for $2\leq p<\frac{p_{0}}{p_{0}-1}
\approx2.18006,$ the following inequalities
\[
\left(  2^{\frac{1}{2}-\frac{1}{p}}\right)  ^{m-1}\leq C_{(m),p}\leq\left(
\frac{1}{\sqrt{2}}\left(  \frac{\Gamma\left(  \frac{2p-1}{2p-2}\right)
}{\sqrt{\pi}}\right)  ^{\frac{1-p}{p}}\right)  ^{m-1}.
\]

For $p\geq2$, we know that
\begin{equation}
\left(  \sum\limits_{j=1}^{\infty}\left(  \sum\limits_{k=1}^{\infty}\left\vert
A(e_{j},e_{k})\right\vert ^{\lambda}\right)  ^{\frac{1}{\lambda}2}\right)
^{\frac{1}{2}}\leq\sqrt{2}\left\Vert A\right\Vert \label{65a}
\end{equation}
with $\lambda=\frac{p}{p-1},$ for all continuous bilinear forms $A:\ell
_{p}\times c_{0}\rightarrow$ $\mathbb{R}$ (see, for instance, \cite[Theorem
1.2 and Remark 5.1]{bayart}). By interpolating (\ref{65a}) and the result of
Theorem \ref{ppp} for $m=2$ in the sense of \cite{bayart} or using the
H\"{o}lder inequality for mixed sums (\cite{bp}) we obtain, for $p\geq
\frac{p_{0}}{p_{0}-1}\approx2.18006,$
\[
\left(  \sum\limits_{j,k=1}^{\infty}\left\vert A(e_{j},e_{k})\right\vert
^{\frac{4p}{3p-2}}\right)  ^{\frac{3p-2}{4p}}\leq\left(  \sqrt{2}\left\Vert
A\right\Vert \right)  ^{1/2}\left(  2^{\frac{p-2}{2p}}\left\Vert A\right\Vert
\right)  ^{1/2}=2^{\frac{1}{2}-\frac{1}{2p}}\left\Vert A\right\Vert .
\]
Using the approach of the previous section we obtain the lower estimate
\[
C_{p,\infty}\geq\frac{\left(  \sum\limits_{j,k=1}^{2}\left\vert T_{2}
(e_{j},e_{k})\right\vert ^{\frac{4p}{3p-2}}\right)  ^{\frac{3p-2}{4p}}
}{\left\Vert T_{2}\right\Vert }=\frac{4^{\frac{3p-2}{4p}}}{2}=2^{\frac{1}
{2}-\frac{1}{p}}
\]
and thus
\[
2^{\frac{1}{2}-\frac{1}{p}}\leq C_{p,\infty}\leq2^{\frac{1}{2}-\frac{1}{2p}}.
\]
When $p=\infty$ we recover the well known optimal estimate of the famous
Littlewood's $4/3$ that can be found in \cite{diniz}.

\medskip

\textbf{Acknowledgement.} The authors are indebted to the two anonymous
referees for their important contributions to the final version of this paper.


\begin{thebibliography}{99}


\bibitem {www}N. Albuquerque, G. Ara\'{u}jo, D. N\'{u}\~{n}ez-Alarc\'{o}n, D.
Pellegrino, and P. Rueda, Bohnenblust--Hille and Hardy--Littlewood
inequalities by blocks, arXiv:1409.6769v6 [math.FA].

\bibitem {bayart}N. Albuquerque, F. Bayart, D. Pellegrino, J.
Seoane-Sepulveda, Sharp generalizations of the multilinear Bohnenblust-Hille
inequality, J. Funct. Anal. \textbf{266} (2014), no. 6, 3726--3740.

\bibitem {a22}N. Albuquerque, F. Bayart, D. Pellegrino and J. B.
Seoane-Sep\'{u}lveda, Optimal Hardy-Littlewood type inequalities for
polynomials and multilinear operators, Isr. J. Math. \textbf{211} (2016), 197--220.

\bibitem {a1}G. Ara\'{u}jo, D.\ Pellegrino, Lower bounds for the constants of
the Hardy-Littlewood inequalities, Linear Algebra Appl. \textbf{463} (2014), 10--15.

\bibitem {bra}G. Ara\'{u}jo, D.\ Pellegrino, On the constants of the
Bohnenblust--Hille and Hardy--Littlewood inequalities, to appear in Bull.
Braz. Math. Soc.

\bibitem {a2}G. Ara\'{u}jo, D.\ Pellegrino, D. Diniz P. Silva e Silva, On the
upper bounds for the constants of the Hardy-Littlewood inequality, J. Funct.
Anal. \textbf{267} (2014), no. 6, 1878--1888.

\bibitem {bohr}F. Bayart, D. Pellegrino and J. B. Seoane-Sep\'{u}lveda, The
Bohr radius of the $n$--dimensional polydisk is equivalent to $\sqrt{(\log
n)/n}$, Adv. Math. \textbf{264} (2014), 726--746.

\bibitem {bp}A. Benedek, R. Panzone, The space $L_{p}$, with mixed norm, Duke
Math. J. \textbf{28 }1961 301--324.

\bibitem {bh}H. F. Bohnenblust, E. Hille, On the absolute convergence of
Dirichlet series, Ann. of Math. \textbf{32} (1931), 600--622.

\bibitem {botelho}G. Botelho, J. Santos, A Pietsch domination theorem for
$(\ell_{p}^{s},\ell_{p})$-summing operators, Arch. Math. (Basel) \textbf{104}
(2015), 47--52.

\bibitem {was}W. Cavalcante, D. N\'{u}\~{n}ez-Alarc\'{o}n, Remarks on the
Hardy--Littlewood inequality for $m$-homogeneous polynomials and $m$-linear
forms, to appear in Quaest. Math.

\bibitem {ann}A. Defant, L. Frerick, J. Ortega-Cerd\'{a}, M. Ouna{\"{\i}}es,
K. Seip, The Bohnenblust-Hille inequality for homogeneous polynomials is
hypercontractive, Ann. of Math. (2), \textbf{174} (2011), 485--497.

\bibitem {Di}J. Diestel, H. Jarchow, A. Tonge, Absolutely summing operators,
Cambridge University Press, Cambridge, 1995.

\bibitem {dimant}V. Dimant and P. Sevilla--Peris, Summation of coefficients of
polynomials on $\ell_{p}$ spaces, Publ. Mat. \textbf{60} (2016), 289--310.

\bibitem {diniz}D. Diniz, G. Mu\~{n}oz-Fern\'{a}ndez, D. Pellegrino, J.
Seoane-Sep\'{u}lveda, Lower bounds for the constants in the Bohnenblust--Hille
inequality: the case of real scalars, Proc. Amer. math. Soc. \textbf{142}
(2014), 575--580.

\bibitem {Ha}U. Haagerup, The best constants in the Khinchine inequality,
Studia Math. \textbf{70} (1982) 231--283.

\bibitem {hardy}G. Hardy and J. E. Littlewood, Bilinear forms bounded in space
$[p,q]$, Quart. J. Math. \textbf{5} (1934), 241--254.

\bibitem {LLL}J. E. Littlewood, On bounded bilinear forms in an infinite
number of variables, Quart. J. Math. \textbf{1} (1930), 164--174.

\bibitem {monta}A. Montanaro, Some applications of hypercontractive
inequalities in quantum information theory. J. Math. Phys. \textbf{53} (2012),
no. 12, 122206, 15 pp.

\bibitem {tonge}B. Osikiewicz, A. Tonge, An interpolation approach to
Hardy--Littlewood inequalities for norms of operators on sequence spaces,
Linear Algebra Appl. \textbf{331} (2001), 1--9.

\bibitem {daniel}D. Pellegrino, The optimal constants of the mixed $(\ell
_{1},\ell_{2})$-Littlewood inequality, J. Number Theory \textbf{160} (2016), 11--18.

\bibitem {popa3}D. Popa, Multiple Rademacher means and their applications, J.
Math. Anal. Appl. \textbf{386} (2012), 699--708.

\bibitem {popa}D. Popa, Multiple summing operators on $\ell_{p}$-spaces,
Studia Math. \textbf{225} (2014), no. 1, 9--28.

\bibitem {rueda}P. Rueda, E.A. S\'{a}nchez-P\'{e}rez, Factorization of
$p$-dominated polynomials through $L_{p}$-spaces, Michigan Math. J.
\textbf{63} (2014), no. 2, 345--353.

\bibitem {velanga}J. Santos and T. Velanga, A note on the Bohnenblust--Hille
inequality for multilinear forms, arXiv:1604.00040v2 [math.FA].

\bibitem {diana}D. M. Serrano-Rodr\'{\i}guez, Improving the closed formula for
subpolynomial constants in the multilinear Bohnenblust--Hille inequalities,
Linear Algebra Appl. \textbf{438} (2013) 3124--3138.
\end{thebibliography}
\end{document}